\documentclass[conference]{IEEEtran}
\usepackage[utf8]{inputenc}
\usepackage{amsmath,amssymb}
\usepackage{graphicx}

\graphicspath{{./},{./images/}}
\newcommand{\equaldist}{\stackrel{d}{=}}

\begin{document}
\title{Congestion analysis of unsignalized intersections}
\author{\IEEEauthorblockN{Abhishek}
\IEEEauthorblockA{Korteweg-de Vries Institute for Mathematics, \\University of Amsterdam, The Netherlands\\
Email: {{\small \tt abhishek@uva.nl}}}
\vspace{2mm}\IEEEauthorblockN{Michel Mandjes}
\IEEEauthorblockA{Korteweg-de Vries Institute for Mathematics, \\University of Amsterdam, The Netherlands\\
Email: {{\small \tt m.r.h.mandjes@uva.nl}}} \and
\IEEEauthorblockN{Marko Boon}
\IEEEauthorblockA{Eindhoven University of Technology,\\ The Netherlands\\
Email: {{\small  \tt m.a.a.boon@tue.nl}}}
\vspace{2mm}
\IEEEauthorblockN{Rudesindo N\'u\~nez-Queija}
\IEEEauthorblockA{Korteweg-de Vries Institute for Mathematics, \\University of Amsterdam, The Netherlands\\
Email: {{\small \tt nunezqueija@uva.nl}}
}
}

\maketitle

\begin{abstract}
This paper considers an unsignalized intersection used by two traffic streams. A stream of cars is using a primary road, and has priority over the other, low-priority, stream. Cars belonging to the latter stream cross the primary road if the gaps between two subsequent cars on the primary road is larger than their critical headways. Questions that naturally arise are: given the arrival pattern of the cars on the primary road, what is the maximum arrival rate of low-priority cars such that the number of such cars remains stable? In the second place, what can be said about the delay experienced by a typical car at the secondary road? This paper addresses such issues by considering a compact model that sheds light on the dynamics of the considered unsignalized intersection. The model, which is of a queueing-theoretic nature, reveals interesting insights into the impact of the user behavior on the above stability and delay issues. The contribution of this paper is twofold. First, we obtain new results for the aforementioned model with driver impatience. Secondly, we reveal some surprising aspects that have remained unobserved in the existing literature so far, many of which are caused by the fact that the capacity of the minor road cannot be expressed in terms of the \emph{mean} gap size; instead more detailed characteristics of the critical headway distribution play a role.
\end{abstract}

\section{Introduction}
A common situation in any road traffic network is that of an unsignalized intersection that is used by two traffic streams which have different priorities. The priority class consists of cars which arrive, according to some inherently random process, at the intersection; the fact that they have priority essentially means that they cross the intersection without observing the low-priority stream. Cars of the low priority stream, however, only cross when the duration (in time) of a gap between two subsequent cars passing by is sufficiently large, i.e., larger than a, possibly car-specific, threshold $T$.

As the high-priority cars do not experience any interference from the low-priority cars, the system's performance is fully determined by the characteristics of the queue of low-priority cars on the secondary road. A first issue concerns the stability of this queue: for what arrival rate of low-priority cars can it be guaranteed that the queue does not explode? Formulated differently: what is the capacity of the minor road? The answer to this question evidently depends on the distribution of the gaps between subsequent cars on the primary road. In addition, however, the specific `preferences' of the low-priority car drivers plays a crucial role: how does the individual car driver choose the threshold $T$ which determines the minimal gap needed. In the existing literature, various models have been studied, the simplest variant being the situation in which all low-priority drivers use the same deterministic $T$. A second, more realistic, model is the one in which the driver draws a random $T$ from a distribution, where the $T$ is resampled for any new attempt; this randomness models the heterogeneity in the preferences of the low-priority car drivers. A further refinement is a model in which the driver sticks to the same $T$ for all his attempts. 

Gap acceptance models are mainly applied to unsignalized intersections (cf. \cite{catchpoleplank,heidemann97,tanner62}), pedestrian crossings (cf. \cite{mayne,tanner51}), and freeways (cf. \cite{drew2,drew1}). Although the gap acceptance process in these three application areas exhibits similar features, the queueing aspects are fundamentally different. In this paper, we focus on motorized vehicles, but all results regarding the capacity of the minor road can easily be applied to pedestrian crossings or freeway merging. Heidemann and Wegmann \cite{heidemann97} give an excellent overview of the existing results in gap acceptance theory, including the three types of user behavior that were discussed; in the sequel we denote these by B$_1$ up to B$_3$. The contribution of the present paper is twofold. First, we include \emph{impatience} of the drivers, waiting to cross the major road, in the model. This phenomenon, which is indeed encountered in practice \cite{abouhenaidy94}, has been studied before (cf. \cite{drew2, drew1, weissmaradudin}), but (to the best of our knowledge) not yet in the context of models B$_2$ and B$_3$, where randomness is encountered in the critical headway $T$. The second contribution of this paper lies in the fact that we reveal some surprising aspects that have remained unobserved in the existing literature so far. In particular we show there is a strict ordering in the capacities, resulting from the different types of driver behavior. This is far from obvious, since it is known that the capacity  {\it not} in terms of the mean quantity ${\mathbb E}[T]$, but more precise distributional information of the random variable $T$ is needed. Perhaps counterintuitively, when comparing two gap time distributions $T_1$ and $T_2$ one could for instance encounter situations in which ${\mathbb E}[T_1]< {\mathbb E}]T_2]$, but in which still the capacity under $T_1$ is smaller than the one under $T_2$.

This paper is structured as follows. In the next section, we describe in more detail the aforementioned types of driver behavior, when waiting for a gap on the main road, and build a model which is used in Section \ref{sect:analysis} to analyze queue lengths and delays. In the subsequent section, we study the traffic congestion on the minor road. In Section \ref{sect:numerics} we present numerical results for several practical examples, focusing on some surprising, paradoxical features that one might encounter. Finally, we give several concluding remarks.

\section{Preliminaries}\label{sect:model}
\begin{figure}[!ht]
\begin{center}
\includegraphics[width=0.7\linewidth]{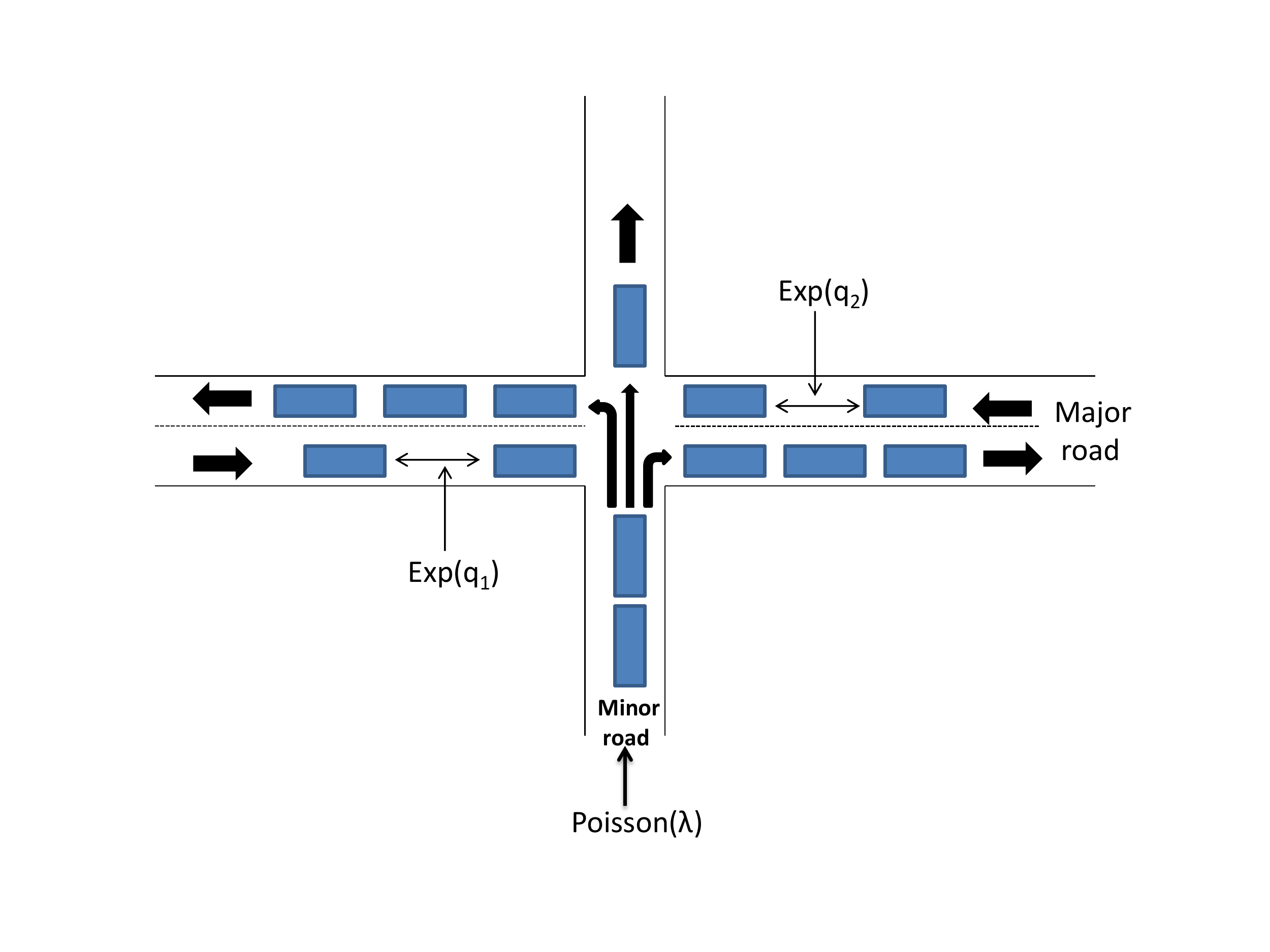}
\end{center}
\caption{The situation analyzed in this paper.}
\label{fig:intersection}
\end{figure}

The situation analyzed in this paper is the following (see also Fig. \ref{fig:intersection}). We consider an intersection used by two traffic streams, both of which wishing to cross the intersection. There are two priorities: the car drivers on the major road have priority over the car drivers on the minor road. The high-priority car drivers arrive at the intersection according to a Poisson process of intensity $q$, meaning that the inter arrival times between any pair of subsequent cars are exponentially distributed with mean $1/q.$ The low-priority car drivers, on the minor road, cross the intersection as soon as they come across a gap with duration larger than $T$ between two subsequent high-priority cars, commonly referred to as the \emph{critical headway}. On the minor road cars arrive according to a Poisson process with rate $\lambda$.

Above we were intentionally imprecise regarding the exact definition of the criterion based on which the low-priority cars decide to cross. Importantly, in this paper we distinguish three types of `behavior' when making this decision.
\begin{itemize}
\item[B$_1$]
The first model is the most simplistic: $T$ is deterministic, and uniform across all low-priority car drivers.
\item[B$_2$]  Clearly B$_1$ lacks realism, in that there will be a substantial level of heterogeneity in terms of driving behavior: one could expect a broad range of `preferences', ranging from very defensive to very reckless drivers. In B$_2$ this is modeled by the car driver at the front end of the queue {\it resampling} $T$ (from a given distribution) at any new attempt (where an `attempt' amounts to comparing this sampled $T$ to the gap between the two subsequent cars that he is currently observing).
\item[B$_3$] In the third model, B$_2$ is made more realistic: it reflects {\em persistent} differences between drivers: each driver selects a random value of $T$, but then sticks to that same value for all attempts, rather than resampling these.
\end{itemize}
For each of these behavior types, we also consider a variant that includes impatience. Details regarding the manner in which impatience is incorporated will be discussed in more detail in the next section.

A few remarks are in place here. In the first place, above we positioned this setup in the context of an unsignalized intersection, but various other applications could be envisioned. One of these could correspond to the situation in which the low-priority cars have to merge with the stream of high-priority cars (e.g. from a ramp). Also in the context of pedestrians crossing a road, the model can be used. We also stress that in the case the primary road actually consists of two lanes that have to be crossed (without a central reservation), with cars arriving (potentially in opposite directions) at Poisson rates (say) $q_1$ and $q_2$, our model applies as well, as an immediate consequence of the fact that  the superposition of two Poisson processes is once again a Poisson process with the parameter $q:=q_1+q_2$.

\section{Analysis of queue lengths and delays}\label{sect:analysis}
The three models can be analyzed using queueing-theoretic techniques. Since results for the variants without impatience have been known in the existing literature (see, for example, Heidemann and Wegmann \cite{heidemann97} for an overview), we will mainly focus on the situation with impatient drivers.

The underlying idea is that we the model can be cast in terms of an M/G/1 queue, i.e., a queue with Poisson arrivals, general service times (which will have their specific form for each of the models B$_1$ up to B$_3$) and a single server. We first introduce some notation. In the first place, we let $X_n$ denote the number of cars in the queue on the minor road when ({\it right after}, that is) the $n$-th low-priority car crosses the primary road; in addition, $T_n^\#$ is the time that this happens.

In our analysis, we need to distinguish between times that the queue  on the secondary road being empty or not. When $X_{n-1}\geqslant 1$ we define the inter departure time between the $(n-1)$-st and $n$-th car from the secondary road as $Y_n:=T_n^\#-T_{n-1}^\#$.
It is clear that for all rules B$_1$ up to B$_3$, we can write $X_n$ as
\[ X_n = \left\{
\begin{array}{l l}
X_{n-1}-1+A(Y_n), & \quad \text{if $X_{n-1} \geqslant 1$ }\\
A(Y_n), & \quad \text{if $X_{n-1} =0$}
\end{array} \right.\]
where $A(\cdot)$ is a Poisson process with intensity $\lambda$.

It is directly verified that, due to the exponentially assumptions imposed, $\{X_n, n=1,2,\ldots\}$ is a Markovian process, the random variables $\{Y_n\}$ are independent identically  distributed (and have a definition for each of the models B$_1$ up to B$_3$). In fact, the process $\{X_n, n=1,2,...\}$ has the dynamics of an M/G/1 queue length process at departures where $Y_n$ corresponds to the $n$-th service time.
As a consequence, we have that $X_n$ has a stationary distribution which is uniquely characterized through its  probability generating function (directly following from  the celebrated Pollaczek--Khinchine formula)
\begin{equation}\label{four}
G_X(z)=\mathbb{E}[z^{X}]=(1-\rho)\frac{(1-z)\tilde{Y}(\lambda(1-z))}{\tilde{Y}(\lambda(1-z))-z}
\end{equation}where $\tilde{Y}(s):=\mathbb{E}[e^{-sY}]$ and  $\rho:=\lambda \mathbb{E}[Y].$ The queue is stable when $\rho < 1$.

Let $W$ and $S$, respectively, be the stationary waiting time and sojourn time corresponding with an arbitrary arriving car on the secondary road.
Since $\{X_n, n=1,2,\ldots\}$  is an M/G/1 queueing length process at departures, the  Laplace-Stieltjes transformations (LSTs) of $W$ and $S$ are given by
\[\tilde{W}(s)=\frac{s(1-\rho)}{\lambda\tilde{Y}(s)+s-\lambda},\:\:\:\tilde{S}(s)=\frac{s(1-\rho)\tilde{Y}(s)}{\lambda\tilde{Y}(s)+s-\lambda},\] respectively.
Now we study the impact of the four types of the driver's behavior on stability and delay.

\vspace{2mm}

\noindent
\textbf{\boldmath B$_1$ {(constant gap):}} Every driver on the minor road needs the same constant critical headway $T_j$ for the $j$-th attempt to enter the main road ($j=1,2,\dots$). We assume that $T_1 \geq T_2\geq \dots \geq T_{\text{min}}$. 
In this case, each user tries a number of attempts (use the memoryless property!),
with success probability ${\mathbb P}(\tau_q>T_j) = e^{-qT_j}$ for the $j$-th attempt, where $\tau_q$ is an
exponential random variable with mean $1/q.$ The Laplace transform of the `service time', ${\mathbb E} [e^{-sY}]$, hence follows from
\begin{equation}
\sum_{k=0}^\infty \left(\prod_{j=1}^k {\mathbb E} [e^{-s\tau_q}1_{\{\tau_q<T_j\}}]\right)\,{\mathbb E} [e^{-sT_{k+1}}1_{\{\tau_q\geqslant T_{k+1}\}}].
\label{lstB1}
\end{equation}
Elementary computations yield
\begin{align*}
{\mathbb E} [e^{-s\tau_q}1_{\{\tau_q<T_j\}}] &= \frac{q}{q+s}\left(1-e^{-(q+s)T_j}\right),\\
{\mathbb E} [e^{-sT_{k+1}}1_{\{\tau_q\geqslant T_{k+1}\}}] &= e^{-(q+s)T_{k+1}},
\end{align*}
which leads to
\begin{equation}
\mathbb{E}[e^{-sY}]=\sum\limits_{k=0}^{\infty}\left(\frac{q}{s+q}\right)^ke^{-(s+q)T_{k+1}}\prod\limits_{j=1}^{k}(1-e^{-(s+q)T_j}).
\label{lstB1final}
\end{equation}
It is also verified that $\mathbb{E}[Y]$ is equal to
\begin{align*}
\sum\limits_{k=0}^{\infty}e^{-qT_{k+1}}\left[ \frac{k}{q}+T_{k+1}-\sum\limits_{i=1}^{k}\frac{T_ie^{-qT_i}}{1-e^{-qT_i}}\right]\prod\limits_{j=1}^{k}(1-e^{-qT_j}).
\end{align*}

These expressions simplify considerably in the case without driver impatience, i.e. $T_1=T_2=\dots=:T$, namely
\begin{align*}
{\mathbb E} [e^{-sY}] &= \frac{(s+q)e^{-(s+q)T}}{s+qe^{-(s+q)T}},\\
{\mathbb E}[Y] &= \frac{e^{qT}-1}{q}.
\end{align*}

\vspace{2mm}

\noindent
\textbf{\boldmath B$_2$ (sampling   per attempt):} With this behavior type, which is also sometimes referred to as ``inconsistent behavior'', every car driver samples a random $T_j$ for its $j$-th `attempt' (where `attempt' corresponds to comparing the resulting $T_j$ with the gap between two subsequent cars on the major road). Although things are slightly more subtle than in B$_1$, with the $T_j$ being random, we can still use the memoryless property of the gaps between successive cars on the major road in combination with the independence of the $T_j$, to argue that expression \eqref{lstB1} is also valid for this model, but now the $T_j$ are random, leading to
\begin{align*}
{\mathbb E} [e^{-s\tau_q}1_{\{\tau_q<T_j\}}] &= \frac{q}{q+s}\left(1-\mathbb{E}[e^{-(q+s)T_j}]\right),\\
{\mathbb E} [e^{-sT_{k+1}}1_{\{\tau_q\geqslant T_{k+1}\}}] &= \mathbb{E}[e^{-(q+s)T_{k+1}}],
\end{align*}
which, after substitution in \eqref{lstB1}, leads to
\[
\mathbb{E}[e^{-sY}]=\sum\limits_{k=0}^{\infty}\left(\frac{q}{s+q}\right)^k\mathbb{E}[e^{-(s+q)T_{k+1}}]\prod\limits_{j=1}^{k}(1-\mathbb{E}[e^{-(s+q)T_j}]).
\]
The mean service time $\mathbb{E}[Y]$ readily follows,
\begin{align*}
&\sum\limits_{k=0}^{\infty}\left[ \mathbb{E}[T_{k+1}e^{-qT_{k+1}}]+\mathbb{E}[e^{-qT_{k+1}}]\left(\frac{k}{q}-\sum\limits_{i=1}^{k}\frac{\mathbb{E}[T_ie^{-qT_i}]}{1-\mathbb{E}[e^{-qT_i}]}\right)\right]\\
&\quad\times \prod\limits_{j=1}^{k}(1-\mathbb{E}[e^{-qT_j}]).
\end{align*}

Again, these expressions simplify considerably in the case without driver impatience.
\begin{align*}
{\mathbb E} [e^{-sY}] &= \frac{(s+q){\mathbb E}[e^{-(s+q)T}]}{s+q\,{\mathbb E}[e^{-(s+q)T}]},\qquad
\mathbb{E}[Y]&=\frac{1-\mathbb{E}[e^{-qT}]}{q\mathbb{E}[e^{-qT}]}. \label{six}
\end{align*}

\noindent
Observe that B$_1$ (with deterministic $T_j$) is a special case of~B$_2$.
\\

\noindent
\textbf{\boldmath B$_3$ (sampling per driver):}   In this variant, sometimes referred to as ``consistent behavior'', every car driver samples a random $T_1$ at his first attempt. This (random) value determines the complete sequence of $T_2, T_3, \dots$. To model this kind of behavior, we introduce a sequence of functions $h_j(\cdots)$, for $j=1, 2, \dots$, such that the critical headway $T_j$ is defined as $T_j:=h_j(T_1)$, where $h_1$ is the identity function. This model combines the advantages of differentiating between various types of drivers, while still modeling a form of impatience.

Since only the first critical headway is random, we obtain the LST of $Y$ by conditioning on the value of $T_1$, and using \eqref{lstB1final}. We obtain
\begin{equation*}
\mathbb{E}[e^{-sY}]=\mathbb{E}\left[\sum\limits_{k=0}^{\infty}\left(\frac{q}{s+q}\right)^ke^{-(s+q)T_{k+1}}\prod\limits_{j=1}^{k}(1-e^{-(s+q)T_j}),\right]
\end{equation*}
where $T_j:=h_j(T_1)$. The mean service time is

\begin{align*}
&\mathbb{E}[Y]=\mathbb{E}\left[
\sum\limits_{k=0}^{\infty}e^{-qT_{k+1}}\left( \frac{k}{q}+T_{k+1}-\sum\limits_{i=1}^{k}\frac{T_ie^{-qT_i}}{1-e^{-qT_i}}\right)\right.\\
&\qquad\qquad\left.\times\prod\limits_{j=1}^{k}(1-e^{-qT_j})\right].
\end{align*}

For completeness, we also give expressions for the variant without driver impatience, which simply result from taking the expectations of the results from model B$_1$.
\begin{align*}
{\mathbb E} [e^{-sY}] &= {\mathbb{E}}\left[\frac{(s+q)e^{-(s+q)T}}{s+qe^{-(s+q)T}}\right],\\
{\mathbb E}[Y] &= \frac{\mathbb{E}[e^{qT}]-1}{q}.
\end{align*}

\vspace{2mm}

\noindent

\section{Impact of driver behavior on capacity}\label{sect:stability}
In the previous section, we have introduced three behavior types, each with and without driver impatience. For these models we determined the distributions of the stationary number of low-priority cars in the queue, being characterized by the corresponding Laplace transform $\tilde{Y}(s)={\mathbb E}[e^{-sY}]$, and hence also the Laplace transform $\tilde{W}(s)$ of the waiting time. In this section we consider the impact of the model of choice (B$_1$ up to B$_3$, that is) on the `capacity' of the secondary road; here `capacity' is defined as the maximum arrival rate $\lambda$ such that the corresponding queue does not explode. A standard result from queueing theory is that for the M/G/1 queue the stability condition is $\rho:=\lambda\,{\mathbb E}[Y]<1$. As a consequence, the capacity of the minor road, denoted by $\bar{\lambda}$, can be determined for each of the models, with or without impatience:
\[
\bar{\lambda}=\frac{1}{{\mathbb E}{[Y]}},
\]
where $\mathbb {E}{[Y]}$ depends on the driver behavior, as explained before. In this section we focus on drivers \emph{without} impatience, denoting the capacity for model B$_i$ by $\bar{\lambda}_i$, for $i=1,2,3$. Although the expressions below can also be found in, for example, Heidemann and Wegmann \cite{heidemann97}, we add some new observations regarding the capacities.

\noindent
\textbf{\boldmath B$_1$ (constant gap):} In this case, we find the stability condition  $\lambda < \bar{\lambda}_1$ where \[\bar{\lambda}_1:=\frac{q}{e^{qT}-1}.\]

\noindent
\textbf{\boldmath B$_2$ (sampling   per attempt):} Using the expression for ${\mathbb E}[Y]$, we now find the condition $\lambda < \bar{\lambda}_2$ where \[\bar{\lambda}_2=\frac{q}{({\mathbb E}[e^{-qT}])^{-1}-1}.\]

\noindent
\textbf{\boldmath B$_3$ (sampling per driver):}
Here, the stability condition is in the form  of  $\lambda < \bar{\lambda}_3$ where \[\bar{\lambda}_3=\frac{q}{\mathbb{E}[e^{qT}]-1}.\]
Importantly, it is here tacitly assumed that the moment generating function $\mathbb{E}[e^{qT}]$ of $T$ exists. A consequence that has not received much attention in the existing literature, is that it also means that in case $T$
has a polynomially decaying tail distribution (i.e., ${\mathbb P}(T>t) \approx C\,t^{-\beta}$ for some $C,\beta>0$ and $t$ large) {\it the queue at the secondary road is never stable}. The reason is that for this type of distributions it is relatively likely that an extremely large $T$ is drawn, such that it takes very long before the car can cross the intersection (such that in the mean time the low-priority queue has built up significantly).

In fact, also for certain light-tailed distributions we find that B$_3$ has an undesirable impact on the capacity. Take, for example, $T$ exponentially distributed with parameter $\alpha$. In this case, we have
\[
\mathbb{E}[Y]=\begin{cases}
1/(\alpha-q), & \quad q < \alpha,\\
\infty & \quad q \geq \alpha,
\end{cases}
\]
implying that the capacity of the minor street drops to zero when $q \geq 1/\mathbb{E}[T]$. Actually, the situation might be even worse than it seems, because it can be shown that  $\mathbb{E}[Y^k]=\infty$ if $q \geq \alpha/k$, for $k=1,2,\dots$. As a consequence, when $\alpha > q \geq \alpha/2$ the capacity is positive, but it follows from \eqref{four} and Little's law that the mean queue length and the mean delay at the minor road grow beyond any bound.\\

Another interesting observation, is that the arrival rates $\bar{\lambda}_1, \bar{\lambda}_2$ and $\bar{\lambda}_3$ obey the ordering
\[\bar\lambda_2\geqslant \bar\lambda_1\geqslant \bar\lambda_3\]
(where in B$_1$ we have chosen $T$ equal to the mean ${\mathbb E}[T]$ used in the other variants). This is an immediate consequence of Jensen's inequality (JI), as we show now. To compare $\bar{\lambda}_1$ and $ \bar{\lambda}_3$ realize that JI implies
\[ \frac{1}{q} (\mathbb{E}[e^{qT}]-1)\geqslant \frac{1}{q}(e^{q\mathbb{E}[T]}-1),\]
which directly entails $\bar{\lambda}_3 \leq \bar{\lambda}_1$. Along the same lines, again appealing to JI,
$
\mathbb{E}[e^{-qT}]\geqslant e^{-q\mathbb{E}[T]},$ and hence
$\bar{\lambda}_2 \geq \bar{\lambda}_1$.

\vspace{2mm}

We conclude this section by stating a number of general observations. In the first place, the above closed-form expressions show that {\em the stability conditions depend on the full distribution of $T$}, as opposed to just the mean values of the random quantities involved.

In many dynamic systems introducing variability degrades the performance of the system. The fact that $\bar\lambda_2 \geqslant \bar\lambda_1$ indicates that in this case this `folk theorem' does not apply: the fact that one resamples $T$ often (every driver selects a new value for each new attempt) actually {\it increases} the capacity of the low-priority road.

\section{Numerical results and practical examples}\label{sect:numerics}

\subsection{Example 1: typical situation}

\begin{figure}[!ht]
\begin{center}
\includegraphics[width=0.7\linewidth]{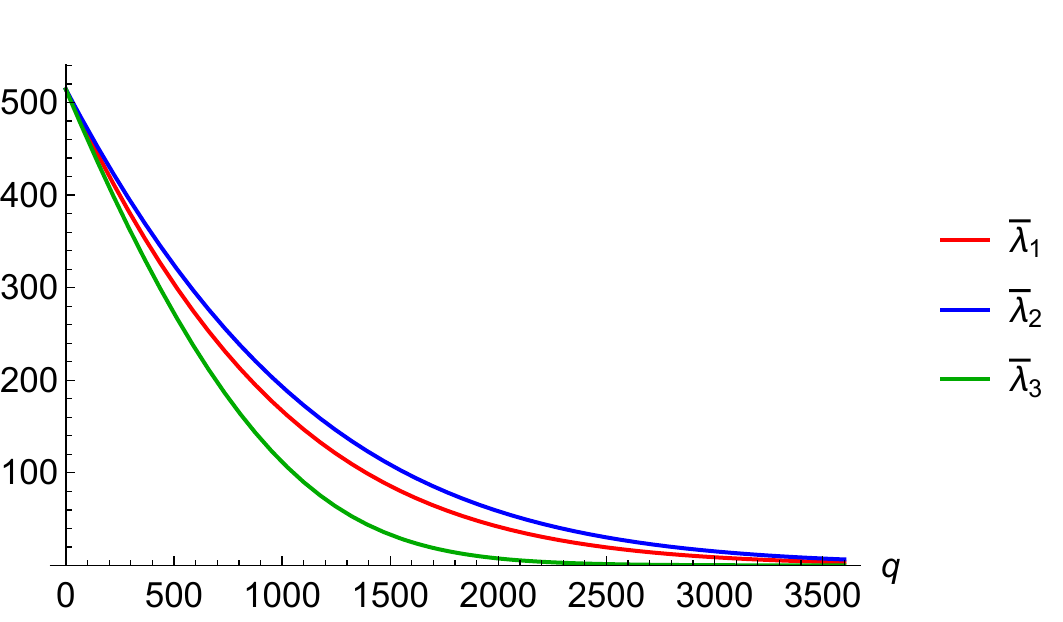}
\end{center}
\caption{Capacity of the minor street (veh/h) as a function of the flow rate on the main road (veh/h) in Example 1.}
\label{fig:example1capacity}
\end{figure}

In this example, in which we try to take realistic parameter settings, we illustrate the impact of driver behavior on the capacity of the system and on the queue lengths. In particular, we compare the following three scenarios (corresponding to the three behavior types):
\begin{itemize}
\item[(1)] All drivers search for a gap between consecutive cars on the major road, that is at least $7$ seconds long.
\item[(2)] A driver on the minor street, waiting for a suitable gap on the major street, will sample a new (random) critical headway every time a car passes on the major street. With probability $9/10$ this critical headway is $6.22$ seconds, and with probability $1/10$ it is exactly $14$ seconds. Note that the expected critical headway is $0.9\times6.22+0.1\times14=7$ seconds, ensuring a fair comparison between this scenario and the previous scenario.
\item[(3)] In this scenario we distinguish between slow and fast traffic. We assume that 90\% of all drivers on the minor road need a gap of (at least) $6.22$ seconds. The other 10\% need at least $14$ seconds.
\end{itemize}
We do not yet incorporate impatience in this example. Fig. \ref{fig:example1capacity} depicts the capacity (veh/h) of the minor street as a function of $q$, the flow rate on the main road (veh/h). The relation $\bar\lambda_2\geq\bar\lambda_1\geq\bar\lambda_3$ is clearly visible. In the next example we will include driver impatience.

\subsection{Example 2: impatience}

We now revisit Example 1, but introducing driver impatience into the model, in the following specific form:
\begin{equation}
T_{k+1}=\alpha (T_k-\Delta)+\Delta, \qquad k = 1, 2, \dots; 0<\alpha<1,\label{impatienceexample}
\end{equation}
which means that the critical headway decreases in every next attempt, approaching the limiting value of $\Delta$. The parameter $\alpha$ determines the speed at which the patience decreases. In scenario 1 all $T_k$ are fixed, with $T_1=7$ seconds. In scenario 2, each of the $T_k$ is a random variable, with $T_1$ equal to $6.22$ or $14$ seconds, with probability $9/10$ and $1/10$ respectively. The distribution of $T_k$ for $k>1$ can be determined from  \eqref{impatienceexample}. Note that the impatience is a new random sample at each attempt, independent of the value of $T_{k-1}$. Scenario 3, as before, is similar to scenario 2, but each driver samples a random impatience $T_1$ exactly once. The value of $T_1$ (which is again either $6.22$ or $14$ seconds) determines the whole sequence of critical gap times at the subsequent attempts according to \eqref{impatienceexample}.

\begin{figure}[!ht]
\parbox{0.49\linewidth}{\centering
\includegraphics[width=\linewidth]{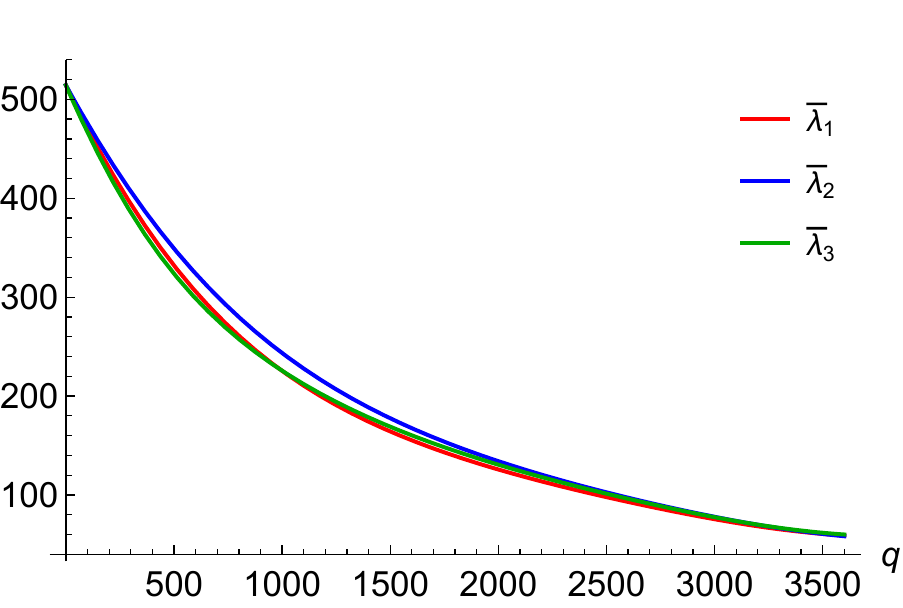}\\
\scriptsize (a) $\alpha=9/10, \Delta=4$ sec.}
\parbox{0.49\linewidth}{\centering
\includegraphics[width=\linewidth]{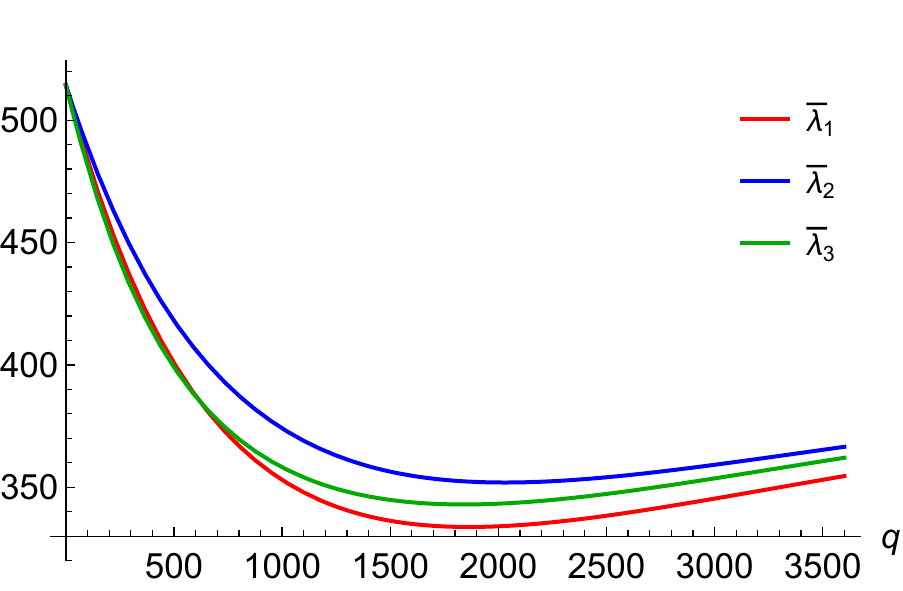}\\
\scriptsize (b) $\alpha=8/10, \Delta=1$ sec.}
\caption{Capacity of the minor street (veh/h) as a function of the flow rate on the main road (veh/h) in Example 2.}
\label{fig:example2capacity}
\end{figure}

Fig. \ref{fig:example2capacity}(a) shows the capacity as a function of $q$, when $\alpha=9/10$ and $\Delta=4$ seconds. It is noteworthy that the strict ordering that was observed in the case without impatience, is no longer preserved, even though $\mathbb{E}[T_k]$ is the same in all scenarios, for $k$ fixed. Another interesting phenomenon, depicted in Fig. \ref{fig:example2capacity}(b), occurs when we decrease the parameter values to $\alpha=8/10$ and $\Delta=1$ second. Now, the capacity actually \emph{increases} when $q$ exceeds a certain threshold. Due to the increase in $q$, gaps between cars on the major road will be smaller in general, but apparently the benefit of having a (much) lower critical headway at each attempt outweighs the disadvantage of having smaller gaps. Similar paradoxical behavior is studied in more detail in the next two examples.

\subsection{Example 3: paradoxical behavior}

\begin{figure}[!ht]
\parbox{0.48\linewidth}{
\centering
\includegraphics[width=\linewidth]{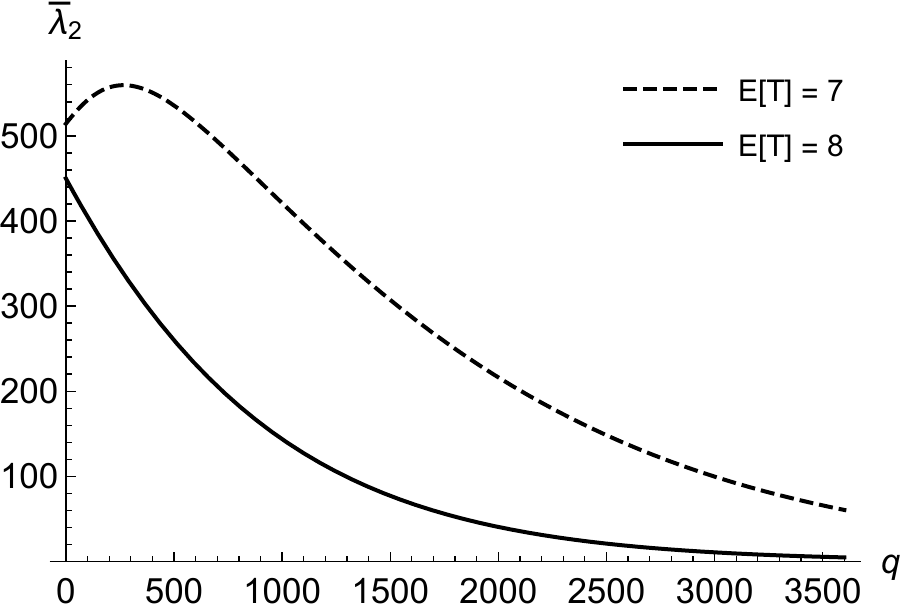}\\
\scriptsize (a) Capacities for model B$_2$.
}
\parbox{0.48\linewidth}{
\centering
\includegraphics[width=\linewidth]{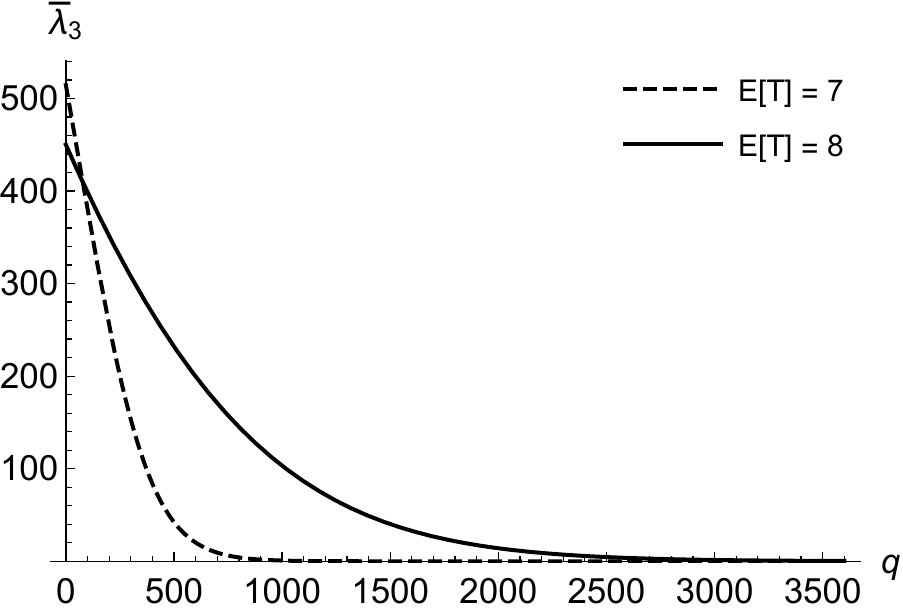}\\
\scriptsize (b) Capacities for model B$_3$.
}
\caption{Capacities $\bar{\lambda}_2$ and $\bar{\lambda}_3$ as a function of $q$ (veh/h) in Example 2.}
\label{fig:example3capacity}
\end{figure}

In this example we illustrate that the gap acceptance model can exhibit unexpected behavior under specific circumstances, caused by the randomness in the behavior of the drivers. To this end, we compare the following two distributions for $T$:
\begin{enumerate}
\item $T_A$ is equal to 4 seconds with probability $9/10$, and 34 seconds with probability $1/10$. The expected critical headway is $\mathbb{E}[T]=7$ seconds.
\item $T_B$ is equal to respectively 6 or 10 seconds, each with probability $1/2$. Now $\mathbb{E}[T]=8$ seconds.
\end{enumerate}
It is apparent that $\mathbb{E}[T_A]<\mathbb{E}[T_B]$, whereas the standard deviation of the first model (9 seconds) is much \emph{greater} than in the second model (2 seconds). In Fig. \ref{fig:example3capacity} we plot the capacities $\bar{\lambda}_2$ and $\bar{\lambda}_3$, when $T$ has the same distributions as $T_A$ and, subsequently, when $T$ has the same distributions as $T_B$. We notice two things. In Fig. \ref{fig:example3capacity}(a), corresponding to model B$_2$, we see that the capacity when $T\equaldist T_A$ first \emph{increases} when $q$ grows larger; only for larger values of $q$ it starts to decrease. This counterintuitive behavior, which is not uncommon in model B$_2$, is explained in more detail in the next numerical example. Now, we mainly focus on the second striking result: in model B$_3$, for $q>78$ veh/h, the capacity when $T\equaldist T_A$ is \emph{less} than the capacity when $T\equaldist T_B$, even though $\mathbb{E}[T_A]<\mathbb{E}[T_B]$. This is a consequence of the fact that the stability condition does not depend on the {\it mean} quantities only, as is common in most queueing models, but it depends on the entire distribution. In this case, the large variance of $T_A$ has a negative impact on the capacity of the minor street.

\subsection{Example 4: The impact of resampling}

\begin{figure}[!ht]
\includegraphics[width=\linewidth]{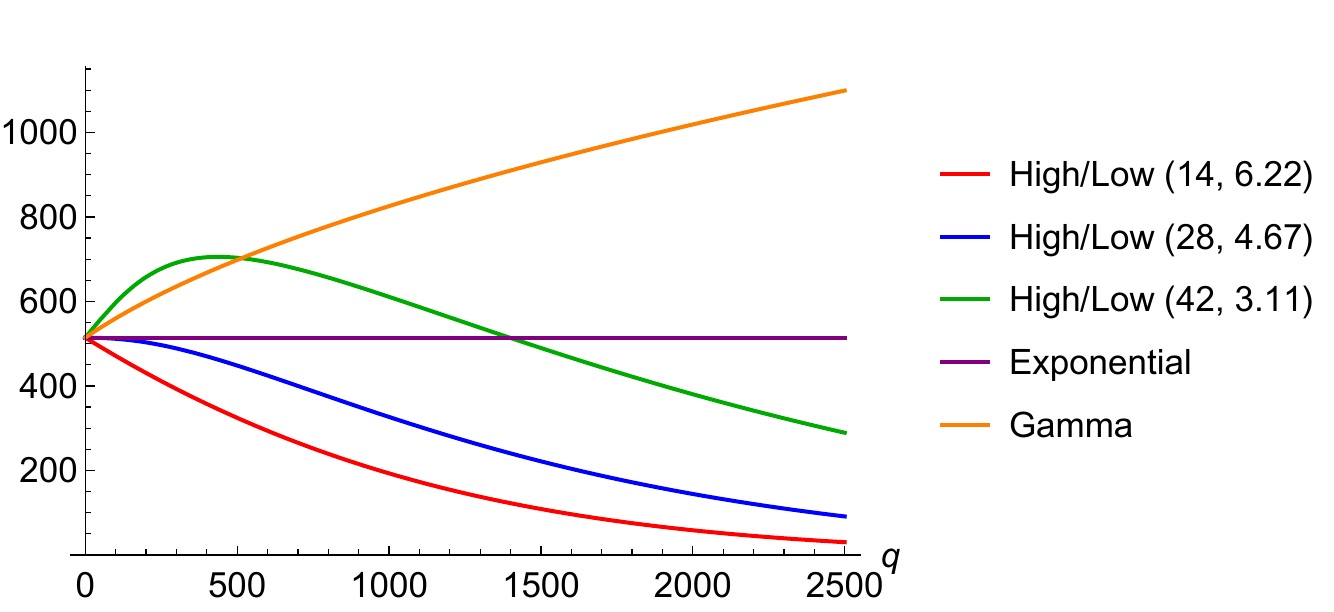}
\caption{Capacity of the minor road (veh/h) for several distributions of $T$.}
\label{fig:example4capacity}
\end{figure}

The previous examples have illustrated that resampling, as described in B$_2$, has a positive impact on the capacity of the minor street. In this example we show that, under specific circumstances, this positive impact may be even bigger than expected. We show this by varying the probability distribution of the critical headway $T$, taking the following five distributions (all with $\mathbb{E}[T]=7$ seconds):
\begin{enumerate}
\item $T$ is equal to 14 seconds with probability $1/10$, and 6.22 seconds with probability $9/10$. This distribution, referred to as High/Low $(14, 6.22)$ in Fig. \ref{fig:example4capacity}, is the same as in Example 1.
\item $T$ is equal to 28 with probability $1/10$, and 4.67 with probability $9/10$. This is a similar distribution as the previous, but with more extreme values.
\item $T$ is equal to 42 with probability $1/10$, and 3.11 with probability $9/10$. Even more extreme values, as in Example 3.
\item $T$ is exponentially distributed with parameter $1/7$.
\item $T$ has a gamma distribution with shape parameter $1/2$ and rate $1/14$.
\end{enumerate}

In Fig. \ref{fig:example4capacity} we have plotted the capacity of the minor road as a function of $q$, the traffic intensity of the main road. Several conclusions can be drawn. First, we notice that in all High/Low distributions, the capacity drops to zero as $q$ increases towards infinity. However, when zooming in at $q$ close to zero, it turns out that model 1 has an immediate capacity drop, whereas model 3 actually \emph{increases} in capacity up to $q\approx 437$. A possible explanation for this striking phenomenon, that we also encountered in the previous numerical example, is that in the third model the critical headway is either extremely low, or extremely high. When a driver has a low critical capacity (which is quite likely to happen), he experiences no delay before leaving the minor street anyway. On the contrast, when a driver has a high critical capacity, he will have to resample, meaning that he has to wait for the next car to pass on the major street. If $q$ increases, the frequency of resampling increases, meaning that he has to wait less before getting a new chance to obtain a small critical headway.

Another extreme case is the exponential distribution, which, due to its memoryless property, has a \emph{constant} capacity, not depending on the traffic flow on the major street.

Clearly, the most paradoxical case is the Gamma distribution with a shape parameter less than one. If $T$ has this particular distribution, it can be shown that the capacity of the minor road keeps on \emph{increasing} as $q$ increases. Although this case, admittedly, may not be a realistic one, it certainly provides valuable insight in the system behavior.

\section{Concluding remarks}\label{sect:conclusions}

Our main target in this work has been to investigate the impact of randomness in the critical headway on the capacity for traffic flows of low priority at a road intersection. For that, we have analyzed three versions of a queueing-theoretic model; each with its own dynamics. Special attention was paid to drivers' impatience under congested circumstances: the value of the critical headway decreases with subsequent attempts to cross the main road.

In our first model (B$_1$) we have assumed that the sequence of critical headways is a deterministically decreasing sequence, and that all cars use same sequence. In the second model (B$_2$), we let each car sample new values for the critical headway, according to a {\em stochastically decreasing} sequence. In the third model (B$_3$) we sample the first value for the critical headway for each car, but then use a deterministic decreasing sequence throughout the attempts of a particular car.

Several instances of our models have previously been studied in literature, but not with the focus on assessing the impact of the drivers' behaviors on the capacity for the low priority floe. Our main observation is that randomness has a strong impact on the capacity. More specifically, the capacity region depends on the entire distribution(s) of the critical headway durations, and not only of the mean value(s). We also observe that resampling of the critical headway values has a benign impact on the capacity of the minor road.

We are currently investigating several extensions to our model. The analysis of our models carries over to variations of it in which, for example, cars may not need the entire duration of their critical headway to cross the main road (the remainder of that duration is the driver's {\em safety margin} to cross the road). Naturally, this further improves the capacity of the minor road, but our main conclusion that capacity is determined by the {\em entire} headway distribution (and not by its mean only) remains equally valid.

\bibliographystyle{abbrv}

\end{document}